\documentclass[12pt]{article}
\usepackage{amsmath}
\usepackage{amsfonts}
\usepackage{amssymb}
\usepackage{amsthm}
\usepackage{stackrel}

\usepackage{color}

\usepackage{figsize}

\usepackage{graphics}

\newtheorem{thm}{Theorem}[section]

\newtheorem{re}{Remark}[section]

\def\binom#1#2{{#1}\choose{#2}}
\newcommand{\R}{{\rm I}\kern-0.18em{\rm R}}
\newcommand{\1}{{\rm 1}\kern-0.25em{\rm I}}
\newcommand{\E}{{\rm I}\kern-0.18em{\rm E}}
\newcommand{\p}{{\rm I}\kern-0.18em{\rm P}}
\makeatletter

\def\@fnsymbol#1{\ensuremath{\ifcase#1\or a\or b\or c\or d\or \e\or f\or *\dagger 	\or \ddagger\ddagger \else\@ctrerr\fi}}

\title{A new convexity-based inequality, characterization of probability distributions and some free-of-distribution tests}
\author{Lev B. Klebanov\footnote{Department of Probability and Mathematical Statistics, Charles University, Prague, Czech Republic. e-mail: levbkl@gmail.com} and Irina V. Volchenkova\footnote{Czech Technical University in Prague}}
\date{}
\begin{document}
\maketitle

\begin{abstract} 

A new inequality between some functional of probability distribution functions is given. The inequality is based on strict convexity of a function used in functional definition. Equality sign in the inequality gives a characteristic property of some probability distributions. This fact together with special character of functional is used to construct free-of-distribution two sample tests.

\end{abstract}

\section{Introduction}\label{sec1} 
\setcounter{equation}{0}

The standing point for this paper is the work by L. Baringhaus and N. Henze \cite{BH}. The first interesting result of theirs work is probabilistic interpretation of the Cram\'{e}r - von Mises distance. Here we provide a generalization of corresponding inequality, which leads to some distribution free tests and to some characterizations of probability distributions. 

\section{Main results}\label{sec2} 
\setcounter{equation}{0}

\begin{thm}\label{th1} 
Let $h$ be strictly convex continuously differentiable function on $[0,1]$ such that $h(0)=0$. Suppose that $F(x)$ and $G(x)$ are continuous probability distribution functions (c.d.f.) on real line $\R^1$. 
Then
\begin{equation}\label{eq1}
\int_{-\infty}^{\infty}h(F(x))dG(x)+\int_{-\infty}^{\infty}h(G(x))dF(x) \geq 2 \int_{0}^{1}h(u)du
\end{equation}
with equality if and only if $F(x)=G(x)$.
\end{thm}

\begin{re}\label{re1}
Condition of convexity of $h(x)$ in Theorem \ref{th1} may be changed by the condition of concavity, but the inequality sign must be chanced by opposite one. This remark remains true for all stuff below.
\end{re}

The right hand side of (\ref{eq1}) allows an probabilistic interpretation. Namely, let $X$ and $Y$ be random variables with p.d.f.'s $F(x)$ and $G(x)$ correspondingly. Under conditions of Theorem \ref{th1} we have
\[ \E h(F(Y)) + \E h(G(X)) \geq 2 \int_{0}^{1}h(u)du \]
with equality if and only if $X \stackrel{d}{=}Y$ that is if and only if $X$ and $Y$ are identically distributed.

Let us note that representation of Cram\'{e}r - von Mises distance obtained in \cite{BH} is a particular case of (\ref{eq1}) for the case $h(u)=u^2$. However, similar interpretation may be provided for more general cases. 
\begin{enumerate}
\item Let $m \geq 2$ be an integer number. Consider $h(u)=u^m$. Suppose that $X_1, \ldots ,X_m$ are independent identically distributed (i.i.d.) as $X$ random variables, and $Y_1, \ldots , Y_m$ are i.i.d. as $Y$ random variables. It is easy to see that
\[ \int_{-\infty}^{\infty}F^m(x)dG(x) + \int_{-\infty}^{\infty}G^m(x)dF(x)=\]
\[ =\p\{\max_{j=1, \ldots ,m}X_j <Y\} +\p\{\max_{j=1, \ldots ,m}Y_j <X\}, \]
what gives probabilistic interpretation of the left hand side of (\ref{eq1}). 

\item Previous case gives us a way to consider polynomial function $h(u)$ of degree $m>1$ having non-negative coefficients. Really, if $h(u)=\sum_{k=1}^{m}c_k u^k$, where $c_k \geq 0$, for all $k=1, \ldots , m$, $c_m>0$. Then
\[ \int_{-\infty}^{\infty}h(F(x))dG(x) + \int_{-\infty}^{\infty}h(G(x))dF(x)=\]
\[ = \sum_{k=1}^{m}c_k \Bigl(\p\{\max_{j=1, \ldots ,k}X_j <Y\} +\p\{\max_{j=1, \ldots ,k}Y_j <X\} \Bigr). \]

\item  Let $h$ be non-negative strictly convex continuously differentiable function on $[0,1]$ such that $h(0)=0$. Suppose that random variables $X_j$ and $Y_j$ are as in item 1. Consider Bernstein polynomial of degree $m$ for the function $h$:
\[ B_m(u) = \sum_{k=1}^{m}\,h(\frac{k}{m})\,{\binom m k}\, u^k\, (1-u)^{m-k}.  \]
It is well known that $B_m$ converges to $h$ uniformly over $[0,1]$ as $m \to \infty$. For statistical interpretation of the inequality (\ref{eq1}) we have to
change $h(u)$ by $B_m(u)$ and pass to limit as $m \to \infty$. Powers of $u$ must be changed by corresponding maximums, while powers of $1-u$ by minimums.
\end{enumerate}

\vspace{3mm}
The result of Theorem \ref{th1} allows a generalization on the case of infinite number of p.d.f.'{s}.

\begin{thm}\label{th2}
Let $h$ be strictly convex continuously differentiable function on $[0,1]$ such that $h(0)=0$. Suppose that $\{F_j(x),\; j=1, 2,\ldots \}$ is a sequence of continuous p.d.f.'s, and $\{p_j,\; j=1, 2, \ldots\}$ is a sequence of positive constants such that $\sum_{j=1}^{\infty}p_j = 1$. Then
\begin{equation}\label{eq2}
\sum_{j\neq k}p_j\cdot p_k \cdot \int_{-\infty}^{\infty}h(F_j(x))dF_k(x) \geq \bigl(1-\sum_{k=1}^{\infty}p_k^2\bigr)\cdot \int_{0}^{1}h(u)du
\end{equation}
with equality if and only if $F_1(x)= \ldots =F_k(x)= \ldots$.
\end{thm}

Inequality (\ref{eq2}) may be generalized on non-countable family of functions $\{F_j\}$. We suppose the reader can make this by her/him self. 

Probabilistic interpretation of the inequality (\ref{eq2}) may be given in the same way as for the case of (\ref{eq1}) and therefore is omitted.

Let us make one trivial remark. Suppose that $Q(x)$ is continuous distribution function. Then (under conditions of Theorem \ref{th2})
\[\sum_{j\neq k}p_j\cdot p_k \cdot \int_{-\infty}^{\infty}h(F_j(x))dF_k(x) = \sum_{j\neq k}p_j\cdot p_k \cdot \int_{-\infty}^{\infty}h(F_j(Q(x)))dF_k(Q(x)).\]
In other words, the value of right hand side of (\ref{eq2}) is invariant with respect to monotone continuous transformation of distribution functions in it. Suppose that we would like to construct a statistical test for the hypothesis $H_o: \; F_1=F_2 = \ldots=F$. Basing on random samples $X_1^{(j)}, \ldots ,X_n^{(j)}$ from population with p.d.f. $F_j$ ($j=1,2, \ldots$) we construct test statistic 
\[T_n=\sum_{j\neq k}p_j\cdot p_k \cdot \int_{-\infty}^{\infty}h(F_j^{(n)}(x))dF_k^{(n)}(x) - \bigl(1-\sum_{k=1}^{\infty}p_k^2\bigr)\cdot \int_{0}^{1}h(u)du, \]    
where $F_j^{(n)}$ is corresponding empirical distribution function. If the hypothesis $H_o$ is true then the distribution of test statistic does not depend on underlying p.d.f. $F$. If $n$ is not too large, we can use computer simulation to find the distribution of the test statistic $T_n$ under $H_o$.

\section{Proof of Theorem \ref{th2}}
\setcounter{equation}{0}

\begin{proof}
It is obvious that Theorem \ref{th1} is a particulat case of Theorem \ref{th2} and, therefore, it is sufficient to give a proof of the second result.

We suppose that $h$ is a strictly convex function. Therefore, if the numbers $p_1, \ldots ,p_n$ ($n \geq 2$) satisfy conditions $p_j>0$ and $\sum_{j=1}^{n}p_j =1$ then 
\begin{equation}\label{eq3}
\sum_{j=1}^{n}p_j h(u_j) \geq h(\sum_{j=1}^{n}p_j u_j)
\end{equation}
for all $u_j \in [0,1]$. The equality in (\ref{eq3}) holds if and only if $u_1 = \ldots = u_n$. Let us take $u_j =F_j(x)$ and integrate both sizes of (\ref{eq3}) over $x$ with respect of p.d.f. $\sum_{j=1}^{n}p_j F_j(x)$. We obtain
\begin{equation}\label{eq4}
\sum_{j=1}^{n}\sum_{k=1}^{n}p_j p_k \int_{-\infty}^{\infty}h(F_j(x))dh(F_k(x)) \geq \int_{0}^{1}h(u)du.
\end{equation} 
However, 
\[ \int_{-\infty}^{\infty}h(F_k(x))dh(F_k(x)) = \int_{0}^{1}h(u)du, \quad k=1, \ldots ,n . \]
Therefore, (\ref{eq4}) implies
\begin{equation}\label{eq5}
\sum_{j\neq k}p_j\cdot p_k \cdot \int_{-\infty}^{\infty}h(F_j(x))dF_k(x) \geq \bigl(1-\sum_{k=1}^{n}p_k^2\bigr)\cdot \int_{0}^{1}h(u)du.
\end{equation}
The inequality (\ref{eq2}) may be obtained as a limit case of (\ref{eq5}) as $n \to \infty$.

From strictly convexity property of $h$ it follows that the equality sign in (\ref{eq2}) holds for equal functions $F_1, \ldots ,F_n, \ldots$ only.
\end{proof}

It is obvious that the case of strictly concave function $h$ is similar, but the inequality sign is opposite. 

\section{One similar inequality}
\setcounter{equation}{0}

Here we propose to change the class of strictly convex functions $h$ by that of positive logarithmically strictly convex functions $\xi$ on $[0,1]$ interval\footnote{That is by class of positive functions haveng strictly convex logarithm}. The function $\xi$ is logarithmically strictly convex if 
\[ \xi((x_1+x_2)/2) \leq \sqrt{\xi(x_1)\xi(x_2)}  \] 
with equality if and only if $x_1=x_2$. This inequality may be rewritten in the form
\begin{equation}\label{eq{6}}
\xi^2((x_1+x_2)/2) \leq \xi(x_1)\xi(x_2)
\end{equation}
From here we get for arbitrary p.d.f.'s $F(x)$ and $G(x)$
\[\int_{-\infty}^{\infty}\xi^2\Bigl(\frac{F(x)+G(x)}{2}\Bigr)d\frac{F(x)+G(x)}{2} \leq \int_{-\infty}^{\infty}\xi(F(x))\xi(G(x))d\frac{F(x)+G(x)}{2}. \]

After simple calculations we obtain
\begin{equation}\label{eq7}
\int_{-\infty}^{\infty} \xi(G(x))d\Xi(F(x)) + \int_{-\infty}^{\infty} \xi(F(x))d\Xi(G(x)) \geq 2 \int_{0}^{1}\xi^2(u)du,
\end{equation}

where $\Xi(u)=\int_{0}^{u}\xi(v)dv$. The equality in (\ref{eq7}) is attend for the case $F(x)=G(x)$ only. Of course, (\ref{eq7}) is very similar to (\ref{eq1}). 

Suppose that $X_1, \ldots ,X_n$ are i.i.d. observations with p.d.f. $F(x)$, and $Y_1, \ldots ,Y_n$ are i.i.d. observations with p.d.f. $G(x)$. Let us consider statistic 
\[\tau_n = \int_{-\infty}^{\infty} \xi(G_n(x))d\Xi(F_n(x)) + \int_{-\infty}^{\infty} \xi(F_n(x))d\Xi(G_n(x)) - \int_{0}^{1}\xi^2(u)du. \]
Here $F_n$ and $G_n$ are corresponding empirical distributions. It is possible to use $\tau_n$ to test the hypothesis $F=G$. It is clear that the distribution of $\tau$ under hypothesis does not depend on the under4lying distribution $F$.

As previously, opposite inequality holds for logarithmically strictly concave function $\xi$.

\section{Acknowledgment}
The work was partially supported by Grant GACR 16-03708S.
Authors are grateful to Dr. K. Helisov\'{a} for her useful comments and kind attention to this work.

\end{document}